\documentclass[12pt,a4paper]{amsart}%
\usepackage{pgf,tikz}
\usepackage{amssymb}
\usepackage{amsmath}
\usepackage{amsfonts}
\usepackage{graphicx}
\usepackage{enumitem}
\usepackage{bm}
\usepackage{comment}
\usepackage{pgfplots}
\usepackage{xurl}
\usepackage{float}
\usepackage[utf8]{inputenc}%
\providecommand{\U}[1]{\protect\rule{.1in}{.1in}}
\usetikzlibrary{backgrounds}
\usetikzlibrary{arrows}
\pgfplotsset{compat=1.15}
\let\oldmathbf\mathbf
\renewcommand{\mathbf}[1]{\boldsymbol{\oldmathbf{#1}}}
\newtheorem{theorem}{Theorem}
\newtheorem*{theorem*}{Theorem}

\newtheorem{definition}[theorem]{Definition}

\newtheorem{lemma}[theorem]{Lemma}

\newtheorem{remark}[theorem]{Remark}

\newcommand{\mycomment}[1]{}
\allowdisplaybreaks
\definecolor{zzttqq}{rgb}{0.6,0.2,0.}
\allowdisplaybreaks
\begin{document}
\title[Discrepancy and absolutely continuous measures]{\textbf{Discrepancy and approximation of absolutely continuous measures with
atomic measures}}
\author[L. Brandolini]{Luca Brandolini}
\address{Dipartimento di Ingegneria Gestionale, dell'Informazione e della Produzione,
Universit\`{a} degli Studi di Bergamo, Viale Marconi 5, 24044 Dalmine BG, Italy}
\email{luca.brandolini@unibg.it}
\author[L. Colzani]{Leonardo Colzani}
\address{Dipartimento di Matematica e Applicazioni, Universit\`{a} di Milano-Bicocca,
Via Cozzi 55, 20125 Milano, Italy}
\email{leonardo.colzani@unimib.it}
\author[G. Travaglini]{Giancarlo Travaglini}
\address{Dipartimento di Matematica e Applicazioni, Universit\`{a} di Milano-Bicocca,
Via Cozzi 55, 20125 Milano, Italy}
\email{giancarlo.travaglini@unimib.it}
\subjclass[2010]{Primary 11K38, 42B10}
\keywords{Geometric discrepancy, Roth's theorem, Fourier transforms, Cassels-Montgomery
lemma, Morrey spaces, Jittered sampling, Gates of Hell Problem}
\maketitle

\begin{abstract}
We prove several results concerning the discrepancy, tested on balls in the
$d$-dimensional torus $\mathbb{T}^{d}$, between absolutely continuous measures
and finite atomic measures.

\end{abstract}

\section{Introduction}

\textit{Is there any measure which is more difficult to approximate by finite
atomic measures with equal weights than the Lebesgue measure?}

\bigskip

This intriguing question has been raised by C. Aistleitner, D. Bilyk, and A.
Nikolov in \cite{ABD} (where it has been termed \textit{Gates of Hell} problem, 
by referring to Rodin's sculpture in Stanford), and of course it can be
specified in more than one way. See also \cite{BCT23,chen85,chen87,FGW,Weiss}.
In this paper we focus on the  related problem of estimating the irregularities
of distribution in the case of absolutely continuous measures on
$\mathbb{T}^{d}$, and we prove sharp bounds for the discrepancies, tested on
balls. Each one of these bounds depends explicitly on the number of points, on
the square norm of the finite sequence of weights associated to the points,
and on the $L^{p}$ or Morrey norm of the function appearing in the absolutely
continuous measure. 

\medskip

The need for an understandable introduction leads us to take a step back.

\medskip

In his celebrated 1954 Mathematika paper \cite{roth} K. Roth proved the
following result, which is the starting point for the field termed
\textit{Geometric Discrepancy}. \medskip

Here and thereafter we identify $\mathbb{T}^{d}=\mathbb{R}^{d} /
\mathbb{Z}^{d}$ with the unit cube $\left[  -1/2,1/2 \right)  ^{d}$.

\begin{theorem}
[Roth]\label{Roththeorem} Let $d\geq2$. There exists a constant $c=c\left(
d\right)  >0$ such that, for every choice of $N$ points $z_{1},\ldots z_{N}$,
not necessarily distinct, in the torus $\mathbb{T}^{d}$, we have the following
lower bound for the $L^{2}$ norm of the discrepancy.
\begin{equation}
\int_{\mathbb{T}^{d}}\left\vert \frac{1}{N}\sum_{j=1}^{N}\chi_{I_{t}}\left(
z_{j}\right)  -\left\vert I_{t}\right\vert \right\vert ^{2}\ dt_{1}\cdots
dt_{d}\geqslant c\,\frac{\log^{d-1}\left(  N\right)  }{N^{2}}\ ,
\label{roth54}%
\end{equation}
where $t=\left(  t_{1},\ldots,t_{d}\right)  \in\left[  0,1\right)  ^{d}$,
$I_{t}=\left[  0,t_{1}\right]  \times\cdots\times\left[  0,t_{d}\right]  $,
$\left\vert G\right\vert $ denotes the volume of a measurable set $G$, and
\begin{equation}
\chi_{G}\left(  x\right)  =\left\{
\begin{array}
[c]{ccc}%
1 &  & \text{if }x\in G,\\
0 &  & \text{if }x\notin G,
\end{array}
\right.  \label{not charact}%
\end{equation}
denotes its characteristic function.
\end{theorem}

%\medskip

H. Davenport \cite{davenport} and K. Roth \cite{roth80} have shown that the
lower bound in (\ref{roth54}) cannot be improved (see also \cite{chen15}).

\medskip

The monographs \cite{BC,matousek} are the basic references for Roth's theorem
and its developments. See also \cite{Chaz00,chenLN,DT,montgomery,travaglini}.

\medskip

Roth's theorem can be rewritten in the following equivalent form (see e.g.
\cite{drmota1,montgomery, ruzsa,schmidt77}), where the parallepipeds anchored
at the origin are replaced by translated and dilated copies of a given cube
$Q$.

\begin{theorem}
Let $d\geq2$ and let $Q=\left[  0,1\right)  ^{d}$. Then, for every choice of
$N$ points $z_{1},\ldots z_{N}$, not necessarily distinct, in the torus
$\mathbb{T}^{d}$, we have%
\[
\int_{0}^{1}\int_{\mathbb{T}^{d}}\left\vert \frac{1}{N}\sum_{j=1}^{N}%
\chi_{\lambda Q+x}\left(  z_{j}\right)  -\lambda^{d}\right\vert ^{2}%
~dxd\lambda\geqslant c \ \frac{\log^{d-1}\left(  N\right)  }{N^{2}}\;.
\]

\end{theorem}

\medskip

A basic development of geometric discrepancy came from a seminal paper by W.
Schmidt \cite{schmidt69}, who observed that replacing parallelepipeds or cubes
with balls leads to a much higher lower bound. This has suggested an
investigation of the relation between the quality of the sampling points
$z_{1},z_{2},\ldots,z_{N}$ and the geometry of the sets (cubes, balls, convex
bodies, fractals, ...) we use to test them. See e.g.
\cite{AMT,BCGGT,BCT23,BT20,BT22,CT,G 2020,matousek}.

H. Montgomery \cite{montgomery} and J. Beck \cite{beck} have sharpened
Schmidt's result and have independently proved the first part of the following
theorem (see also \cite{BCT97}). See e.g. \cite{BGT,BIT,Ken} for proofs of the
second part.

\begin{theorem}
[Schmidt - Montgomery - Beck]\label{SMB}Let $d\geq2$ and let $B\subset
\mathbb{T}^{d}$ be a ball of diameter $<1$. Then there exists a constant
$c_{1} >0$ such that, for every choice of $N$\ points $z_{1},\ldots,z_{N}$,
not necessarily distinct, in $\mathbb{T}^{d}$, we have%
\[
\int_{0}^{1}\int_{\mathbb{T}^{d}}\left\vert \frac{1}{N}\sum_{j=1}^{N}%
\chi_{\lambda B+x}\left(  z_{j}\right)  -\lambda^{d}\left\vert B\right\vert
\right\vert ^{2}~dxd\lambda\geqslant c_{1} N^{-1-1/d}\;,
\]
where $\left\vert B\right\vert $ is the volume of the ball $B$. Moreover, for
every $0<\lambda\leqslant1$ there exists a constant $c_{2}$ such that, for
every positive integer $N$, there exist points $\widetilde{z}_{1}%
,\ldots,\widetilde{z}_{N}$ in $\mathbb{T}^{d}$ such that%
\[
\int_{\mathbb{T}^{d}}\left\vert \frac{1}{N}\sum_{j=1}^{N}\chi_{\lambda
B+x}\left(  \widetilde{z}_{j}\right)  -\left\vert B\right\vert \right\vert
^{2}~dx\leqslant c_{2} N^{-1-1/d}\ .
\]

\end{theorem}

\medskip

Perhaps it is interesting to recall that the first part of Theorem \ref{SMB}
does not hold true if we forget the integration in the dilation variable
$\lambda$, that is we fix $\lambda$. Indeed, in this case the results may
depend on the dimension $d$ in an unexpected way. See e.g.
\cite{BCGT2015,BGT,PaSo}. See also \cite{TrTu} for a curious phenomenon
concerning rotated polygons.

\medskip

As in \cite{chen20}, we reconsider the discrepancy in the following more
general setting.

\medskip

\textit{We want to estimate the quality of the distribution of }%
$N$\textit{\ points }$z_{1},\ldots,z_{N}$,\textit{ not necessarily distinct,
with respect to a Borel probability measure }$\mu$ \textit{on} $\mathbb{T}%
^{d}$. \textit{We choose a reasonably large family $\mathcal{R}$ of measurable
sets and, for $R\in\mathcal{R}$, we study the discrepancy}%
\[
D_{N}\left(  R\right)  =\frac{1}{N}\sum_{j=1}^{N}\chi_{R}\left(  z_{j}\right)
-\mu\left(  R\right)  \;.
\]

In \cite{BCT23} we have investigated the case of possibly singular measures
$\mu$ and a family $\mathcal{R}$ of possibly fractal sets. Here we move in the
opposite direction and consider absolutely continuous measures $f(x)dx$. Hence
we meet again the question raised at the beginning of this paper, which we
modify as follows.

\smallskip

\textit{Is there any absolutely continuous measure $f\left(  x\right)  dx$
which is easier to approximate by finite atomic measures than the Lebesgue
measure $dx$?}

\medskip

In \cite{chen85} W. Chen provided a negative answer.

\begin{theorem}
[Chen]\label{chen}Let $f(x)$ be a non-zero Lebesgue integrable real-valued
function on $\mathbb{T}^{d}$ and let $h(x)$ be any real-valued function on
$\mathbb{T}^{d}$. Let us consider $N$\textit{\ points }$z_{1},\ldots,z_{N}%
$,\textit{ not necessarily distinct in }$[0,1]^{d}$. For every $t=\left(
t_{1},\ldots,t_{d}\right)  \in\left[  0,1\right)  ^{d}$ we consider the
parallelepiped (anchored at the origin)
\[
\mathcal{R}_{t}=\left\{  \left(  u_{1},u_{2},\ldots,u_{d}\right)
:\ 0\leqslant u_{j}\leqslant t_{j},\ \ j=1,2,\ldots,N\right\}  \ .
\]
Then there exists a constant $c>0$ depending only on the dimension $d$ and on
the function $f$ such that
\begin{equation}
\int_{\mathbb{T}^{d}}\left\vert \frac{1}{N}\sum_{j=1}^{N}\chi_{\mathcal{R}%
_{t}}(z_{j})h(z_{j})-\int_{\mathbb{T}^{d}}\chi_{\mathcal{R}_{t}}%
(x)f(x)\ dx\right\vert ^{2}dt\geqslant c \ \frac{\log^{d-1}(N)}{N^{2}}\ .
\label{chen constant}%
\end{equation}

\end{theorem}

Observe that by letting $h(u)=f(u)=\text{constant}\neq0\ $we are back to
Roth's theorem.

Actually Chen's result is more general, since it deals with the $L^{p}$ norm
($p>1$) of the discrepancy.

\medskip

In this paper we are going to prove an analog of Chen's result in the spirit
of the Schmidt-Montgomery-Beck theorem, that is, testing the discrepancy on
balls or on convex bodies having boundaries with positive curvature. We point
out that the constant $c$ in (\ref{chen constant}) depends on $f$ and, in our
setting, one of our goals is to make such a dependence explicit. More
generally, we prove several sharp results on the discrepancies associated to
absolutely continuous measures. See Section \ref{Main} below.

\medskip

Our arguments are essentially Fourier analytic and (a variation of)\ a
classical lemma of Cassels and Montgomery (see \cite[Chapter 6]{montgomery},
see also \cite{BCT23,BBG}) is a basic tool.

\section{Notation}

\begin{itemize}
[leftmargin=2.0cm,labelsep=0.5cm]

\item[$c,c_{1},c_{2},\ldots$ \ ] positive constants which may change from step
to step.

\item[$E\approx F$] when both $E\leqslant c_{1}F$ and $F\leq c_{2}E$ are true.

\item[$\left\vert G\right\vert $] volume of a measurable set $G$.

\item[$\left\Vert \alpha\right\Vert $] normalized $\ell^{2}$ norm of the
weights, see (\ref{not alfanorm}).

\item[$B_{r}$] $=\left\{  x\in\mathbb{T}^{d}:\left\vert x\right\vert
<r\right\}  $, ball of radius $r$.

\item[$B\left(  z,r\right)  $ ] $=\left\{  x\in\mathbb{T}^{d}:\left\vert
x-z\right\vert <r\right\}  $, ball of radius $r$ centered at $z$.

\item[$\left\Vert f\right\Vert _{p}$] $L^{p}$ space norm, see
(\ref{def_norma_Leb}).

\item[$\left\Vert f\right\Vert _{1,\lambda}$] Morrey norm, see (\ref{morrey}).

\item[$\operatorname{diam}\left(  G \right)  $] diameter of a set $G$.

\item[$\mathcal{D}_{N}\left(  x\right)  $] discrepancy (translations), see
(\ref{not discr no dilat}).

\item[$\mathcal{D}_{N}\left(  x,r\right)  $] discrepancy (translations and
dilations), see (\ref{not discr dilat}).

\item[$\mathcal{D}_{N}\left(  x,r,\sigma\right)  $] discrepancy (translations,
dilations, and rotations) see (\ref{not discr dilat rot}).

\item[$\widehat{\mathcal{D}}_{N}\left(  m,r\right)  $] Fourier coefficient of
$\mathcal{D}_{N}\left(  x,r\right)  $, see (\ref{not four}).

\item[$E_{j}$] see (\ref{Ej}).

\item[$\left\Vert f\right\Vert _{p}$] $L^{p}$ space norm, see
(\ref{def_norma_Leb}).

\item[$\left\Vert f\right\Vert _{1,\lambda}$] Morrey norm, see (\ref{morrey}).

\item[$f_{z}$] $=\left(  f\left(  z_{1}\right)  ,\ldots, f\left(
z_{N}\right)  \right)  $, see (\ref{fz}).

\item[$J\left(  N,f,r \right)  $] see (\ref{JNf}).

\item[$K_{M}\left(  x\right)  $] see (\ref{K_M}).

\item[$\mathcal{M}_{E_{j}}$] see (\ref{MQg}).

\item[$\chi_{G}\left(  x\right)  $] characteristic function of a set $G$, see
(\ref{not charact}).
\end{itemize}

\section{Main results}

\label{Main}

We first define the discrepancy that we are going to consider in this paper.

\begin{definition}
\label{defdiscr} Let $f\in L^{1}\left(  \mathbb{T}^{d}\right)  $ and let
$z_{1},\ldots,z_{N}$ be $N$ points in $\mathbb{T}^{d}$, not necessarily
distinct. Let $\alpha_{1},\ldots, \alpha_{N}$ be non-negative real numbers (we
call them weights), let $\alpha=\left(  \alpha_{1},\ldots,\alpha_{N}\right)
$, and let%
\begin{equation}
\left\Vert \alpha\right\Vert =\left\{  \frac{1}{N}\sum_{j=1}^{N}\alpha_{j}%
^{2}\right\}  ^{1/2}\ . \label{not alfanorm}%
\end{equation}
For $x\in\mathbb{T}^{d}$, $0<r<1/2$, and $B_{r}=\left\{  x\in\mathbb{T}%
^{d}:\left\vert x\right\vert <r\right\}  $, we define the discrepancy%
\begin{equation}
\mathcal{D}_{N}\left(  x,r\right)  =\frac{1}{N}\sum_{j=1}^{N}\alpha_{j}%
\chi_{-x+B_{r}}\left(  z_{j}\right)  -\int_{-x+B_{r}}f\left(  y\right)  dy \ .
\label{not discr dilat}%
\end{equation}
When $B$ is a ball with fixed radius, we write $\mathcal{D}_{N}\left(
x\right)  $ in place of $\mathcal{D}_{N}\left(  x,r\right)  $, that is%
\begin{equation}
\mathcal{D}_{N}\left(  x\right)  =\frac{1}{N}\sum_{j=1}^{N}\alpha_{j}\chi
_{B}\left(  z_{j}+x\right)  -\int_{B}f\left(  y+x\right)  dy\ .
\label{not discr no dilat}%
\end{equation}

\end{definition}

\medskip

The following theorem is our first result. Here
\begin{equation}
\left\Vert f \right\Vert _{p} = \left\{  \int_{\mathbb{T}^{d}} \left|
f(t)\right|  ^{p} dt\right\}  ^{1/p}. \label{def_norma_Leb}%
\end{equation}

\begin{theorem}
\label{five}Let $1<p\leqslant2$ and let $1/p+1/q=1$. Let $0<a<b<1/2$. Then
there exists a constant $c>0$ such that for every complex-valued function
$f\in L^{p}\left(  \mathbb{T}^{d}\right)  $, for every choice of $N$ points
$z_{1},z_{2},\ldots,z_{N}\ $in $\mathbb{T}^{d}$, not necessarily distinct, and
non-negative weights $\left\{  \alpha_{j}\right\}  _{j=1}^{N}$, not all zero,
we have%
\begin{equation}
\int_{a}^{b}\int_{\mathbb{T}^{d}}\left\vert \mathcal{D}_{N}\left(  x,r\right)
\right\vert ^{2}dxdr\geqslant c\ N^{-1-q/(2d)}\left\Vert \alpha\right\Vert
^{2+q/d}\left\Vert f\right\Vert _{p}^{-q/d}\ . \label{pq}%
\end{equation}

\end{theorem}

In particular, when $p=2$ we have%
\begin{equation}
\int_{a}^{b}\int_{\mathbb{T}^{d}}\left\vert \mathcal{D}_{N}\left(  x,r\right)
\right\vert ^{2}dxdr\geqslant c\ N^{-1-1/d}\left\Vert \alpha\right\Vert
^{2+2/d}\left\Vert f\right\Vert _{2}^{-2/d}\ . \label{caso P 2}%
\end{equation}

Hence, because of Theorem \ref{SMB}, the above result can be considered a contribution to the question
raised at the very beginning of this paper.

\bigskip

The next result shows that the RHS in (\ref{pq}) cannot be improved if $p=2$.
We do not know if (\ref{pq}) is sharp in the range $1<p<2$. Note that if
$p\rightarrow1^{+}$ and $\left\Vert \alpha\right\Vert =\left\Vert f\right\Vert
_{p}=1$, then the RHS in (\ref{pq}) vanishes. Actually, the second part of the
next theorem shows that for $p=1$ the discrepancy can be as small as we wish.

\begin{theorem}
\label{Lpalto} 1) Let $p>1$. Let $H$ be a positive integer, let $N=H^{d}$ and
let%
\begin{equation}
\left\{  z_{1},\ldots,z_{N}\right\}  =\frac{1}{H}\mathbb{Z}^{d}\cap\left[
-\frac{1}{2},\frac{1}{2}\right)  ^{d}. \label{grid}%
\end{equation}
Assume $\alpha_{1}=\cdots=\alpha_{N}=1$ (hence $\left\Vert \alpha\right\Vert
_{2}=1$). Then there exist non-negative functions $f\in L^{p}\left(
\mathbb{T}^{d}\right)  $ of arbitrarily large $L^{p}$ norms, such that, for
any given ball $B$ in $\mathbb{T}^{d}$,%
\begin{equation}
\int_{\mathbb{T}^{d}}\left\vert \mathcal{D}_{N}\left(  x\right)  \right\vert
^{2}dx\leqslant c\ N^{-1-1/d}\left\Vert f\right\Vert _{p}^{-q/d}\ .
\label{noconvex}%
\end{equation}
2) Let $p=1$ and let $z_{1},\ldots,z_{N}$ be as in (\ref{grid}). Let $\left\{
\varepsilon\left(  N\right)  \right\}  _{N=1}^{+\infty}$ be any positive
vanishing sequence. Then, for every $N$ there exist non-negative functions
$f\in L^{1}\left(  \mathbb{T}^{d}\right)  $, with $\left\Vert f\right\Vert
_{1}=1$, such that%
\begin{equation}
\int_{\mathbb{T}^{d}}\left\vert \mathcal{D}_{N}\left(  x,r\right)  \right\vert
^{2}dx\leqslant\varepsilon\left(  N\right)  \ . \label{solouno}%
\end{equation}

\end{theorem}

\bigskip

The assumption $\alpha_{j}\geqslant0$ for every $j$ in Theorem \ref{five}
cannot be dropped. See Remark \ref{segni} in the next section.

\begin{remark}
\label{Remark Delta} The exponents in the term $\left\Vert \alpha\right\Vert
^{2+q/d}\left\Vert f\right\Vert _{p}^{-q/d}$ in (\ref{pq}) may have the
following explanation. First we observe that the discrepancy associated to the
Dirac delta $\delta_{0}$ can be zero. Indeed, if $z_{1}=\cdots=z_{N}=0$ and
$N^{-1}\sum_{j=1}^{N}\alpha_{j}=1$, then, for any measurable set $R
\subset\mathbb{T}^{d}$, we have%
\[
D_{N}\left(  R\right)  =\chi_{R}\left(  0\right)  -\delta\left(  R\right)
=0\ .
\]
The above trivial observation suggests that the \textquotedblleft
closer\textquotedblright\ is a function $f$ to a Dirac delta, the smaller its
discrepancy can be. Since an $L^{p}$ function close to the Dirac delta has
large $L^{p}$ norm (if $p>1$), then it is natural to expect that a small
discrepancy should be associated to a negative power of the $L^{p}$ norm of
the function $f$. The exponent $2+q/d$ of $\left\Vert \alpha\right\Vert $
comes from the following simple homogeneity argument. Assume we multiply both
the vector $\alpha$ and the function $f$ by the same constant $t$. Then
\[
\int_{a}^{b}\int_{\mathbb{T}^{d}}\left\vert \mathcal{D}_{N}\left(  x,r\right)
\right\vert ^{2}dxdr
\]
gets multiplied by $t^{2}$. Therefore the sum of the exponents of $\left\Vert
\alpha\right\Vert $ and $\left\Vert f\right\Vert _{p}$ must be $2$. In
particular, if $\left\Vert f\right\Vert _{p}$ is raised to the power $-q/d$,
then $\left\Vert \alpha\right\Vert $ needs the power $2+q/d$.
\end{remark}

\medskip

Besides $L^{p}$ spaces, we can consider other function spaces. We are going to
consider Morrey spaces (see e.g. \cite{Adams}), which seem to be quite
tailored for the problems considered in this paper.

\begin{definition}
Let $B\left(  z,r\right)  $ be the ball centered at $z \in\mathbb{T}^{d}$,
with radius $r$, let $0\leqslant\lambda\leqslant d$, and let $f\in
L^{1}\left(  \mathbb{T}^{d}\right)  $. The Morrey spaces $\mathcal{M}%
_{1,\lambda}\left(  \mathbb{T}^{d}\right)  $ are defined by the norms%
\begin{equation}
\left\Vert f\right\Vert _{1,\lambda}=\sup_{z\in\mathbb{T}^{d}, \ 0<r<1/2}%
\ \frac{1}{r^{\lambda}}\int_{B\left(  z,r\right)  }\left\vert f\left(
x\right)  \right\vert dx\ . \label{morrey}%
\end{equation}

\end{definition}

Observe that $\mathcal{M}_{1,0}\left(  \mathbb{T}^{d}\right)  =L^{1}\left(
\mathbb{T}^{d}\right)  $ and, by the Lebesgue differentiation theorem,
$\mathcal{M}_{1,d}\left(  \mathbb{T}^{d}\right)  =L^{\infty}\left(
\mathbb{T}^{d}\right)  $. Also observe that the H\"{o}lder inequality yields
the imbedding $L^{p}\left(  \mathbb{T}^{d}\right)  \subseteq\mathcal{M}%
_{1,d/q}\left(  \mathbb{T}^{d}\right)  $, where $1<p<+\infty$, and $p$, $q$
are conjugate indices. Indeed,
\begin{align}
\left\Vert f\right\Vert _{1,d/q}  &  =\sup_{z\in\mathbb{T}^{d},\ 0<r<1/2}%
\frac{1}{r^{d/q}}\int_{B\left(  z,r\right)  }\left\vert f\left(  x\right)
\right\vert dx\label{imbedding}\\
&  \leqslant\sup_{z\in\mathbb{T}^{d},\ 0<r<1/2}\frac{1}{r^{d/q}}\left\vert
B\left(  z,r\right)  \right\vert ^{1/q}\left\Vert f\right\Vert _{p}%
=c\left\Vert f\right\Vert _{p}\ .\nonumber
\end{align}

The following result concerns the discrepancies associated to functions in
Morrey spaces.

\begin{theorem}
\label{morreybasso}Let $0<a<b<1/2$ and $0<\lambda\leqslant d$. Then there
exists a constant $c>0$ such that for every choice of $N$ points $z_{1}%
,z_{2},\ldots,z_{N}$, not necessarily distinct, in $\mathbb{T}^{d}$, for all
choices of non-negative weights $\alpha_{1},\ldots,\alpha_{N}$ and for every
non-negative function $f\in\mathcal{M}_{1,\lambda}\left(  \mathbb{T}%
^{d}\right)  $, we have%
\begin{equation}
\int_{a}^{b}\int_{\mathbb{T}^{d}}\left\vert \mathcal{D}_{N}\left(  x,r\right)
\right\vert ^{2}dxdr\geqslant cN^{-1-1/\lambda}\left\Vert \alpha\right\Vert
^{2+1/\lambda}\left\Vert f\right\Vert _{1,\lambda}^{-1/\lambda}\ .
\label{morreysimple}%
\end{equation}

\end{theorem}

\begin{remark}
The above theorem and the imbedding of $L^{p}\left(  \mathbb{T}^{d}\right)  $
in the Morrey space $\mathcal{M}_{1,d/q}\left(  \mathbb{T}^{d}\right)  $ (see
(\ref{imbedding})) yield
\begin{align}
\int_{a}^{b}\int_{\mathbb{T}^{d}}\left\vert \mathcal{D}_{N}\left(  x,r\right)
\right\vert ^{2}dxdr  &  \geqslant cN^{1-q/d}\left\Vert \alpha\right\Vert
^{2+q/d}\left\Vert f\right\Vert _{1,d/q}^{-q/d}\label{qfratttod}\\
&  \geqslant cN^{1-q/d}\left\Vert \alpha\right\Vert ^{2+q/d}\left\Vert
f\right\Vert _{p}^{-q/d}\ .\nonumber
\end{align}
The above argument is quite simple, and it is not surprising that the
inequality (\ref{qfratttod}), with the term $N^{1-q/d}$, is weaker than the
inequality (\ref{pq}), with the term $N^{1-q/\left(  2d\right)  }$.
\end{remark}

\medskip

The next result shows that the lower bound (\ref{morreysimple}) is sharp.

\begin{theorem}
\label{morreyalto}Let $0<\lambda<d$, let $H$ be a positive integer, and let
$N=H^{d}$. Let%
\[
\left\{  z_{1},\ldots,z_{N}\right\}  =\frac{1}{H}\mathbb{Z}^{d}\cap\left[
-\frac{1}{2},\frac{1}{2}\right)  ^{d} \ ,
\]
and let $\alpha_{1}=\cdots=\alpha_{N}=1$ (hence $\left\Vert \alpha\right\Vert
=1$). Then there exist functions $f\in\mathcal{M}_{1,\lambda}\left(
\mathbb{T}^{d}\right)  $ of arbitrarily large norms, such that, for every ball
$B$,
\[
\int_{\mathbb{T}^{d}}\left\vert \mathcal{D}_{N}\left(  x\right)  \right\vert
^{2}dx\leqslant cN^{-1-1/\lambda}\left\Vert f\right\Vert _{1,\lambda
}^{-1/\lambda}\ .
\]

\end{theorem}

\begin{remark}
The above results can be extended with no difficulties to the case where the
balls $B_{r}$ are replaced with the dilates of a convex body $C\subset
\mathbb{T}^{d}$, and smooth boundary with positive curvature. The key
ingredients are classical estimates on the decay of the Fourier transform of a
convex body, see e.g. \cite{Stein} and the references therein. If we drop the
smoothness and curvature assumptions, we can obtain similar results by
averaging over rotations of the body (which then needs to have diameter $<1$).
More precisely, for $\sigma\in SO(d)$ we define%
\begin{equation}
\mathcal{D}_{N}\left(  x,r,\sigma\right)  =\frac{1}{N}\sum_{j=1}^{N}\alpha
_{j}\chi_{r\sigma(C)}\left(  z_{j}+x\right)  -\int_{r\sigma(C)}f\left(
y+x\right)  dy\ . \label{not discr dilat rot}%
\end{equation}
Then (\ref{pq}) can be replaced with%
\[
\int_{SO(d)}\int_{a}^{b}\int_{\mathbb{T}^{d}}\left\vert \mathcal{D}_{N}\left(
x,r,\sigma\right)  \right\vert ^{2}\ dxdrd\sigma\geqslant cN^{-1-q/(2d)}%
\left\Vert \alpha\right\Vert ^{2+q/d}\left\Vert f\right\Vert _{p}^{-q/d}\ ,
\]
and (\ref{morreysimple}) can be replaced with%
\[
\int_{SO(d)}\int_{a}^{b}\int_{\mathbb{T}^{d}}\left\vert \mathcal{D}_{N}\left(
x,r,\sigma\right)  \right\vert ^{2}\ dxdrd\sigma\geqslant cN^{-1-1/\lambda
}\left\Vert \alpha\right\Vert ^{2+1/\lambda}\left\Vert f\right\Vert
_{1,\lambda}^{-1/\lambda}\ .
\]
Moreover, there are analogues of Theorem \ref{Lpalto} and Theorem
\ref{morreyalto}. We omit the details. The key observation is that the
classical estimates for the pointwise decay of the Fourier transforms of a
convex body can be replaced by $L^{2}$ spherical averages of the decay of
these Fourier transforms. See \cite{BCT97,BCT23,Bra-Hof-Ios03}.
\end{remark}

\medskip

For the classical case, that is $f\left(  x\right)  dx$ replaced with the
Lebesgue measure $dx$, sharp estimates for the discrepancy can be obtained by
using the so-called \textit{jittered sampling}, which in this paper consists
in partitioning $\mathbb{T}^{d}$ into $N$ disjoint measurable sets $E_{1} ,
\ldots, E_{N}$, having volume $N^{-1}$ and diameter $\approx N^{-1/d}$ (as an
example, assume $N=H^{d}$, where $H$ is a positive integer, and split
$\mathbb{T}^{d}=\left[  -1/2,1/2\right)  ^{d}$ as the union of $N$ equal
half-open disjoint cubes with side-length $N^{-1/d}$ and volume $N^{-1}$). We
recall that the existence of this decomposition is known in the very general
setting of metric measure spaces (see \cite{GL}). Next we take a random point
$z_{j}$ inside each $E_{j}$, and we let $\alpha_{j} = f\left(  z_{j}\right)
$. See e.g. \cite{BCCGT,BGT,ChTr,Doerr,KiPa,PaSt,travaglini}.

\medskip

We have the following result.

\begin{theorem}
\label{jitt} Let $N$ be a positive integer, and consider a decomposition of
$\mathbb{T}^{d}$ into a disjoint union of measurable sets $E_{j}$ satisfying
\begin{equation}
\mathbb{T}^{d}=\bigcup_{j=1}^{N} E_{j} \ , \ \ \text{with} \ \left\vert E_{j}
\right\vert =\frac{1}{N} \ \ \text{and} \ \ \operatorname{diam}\left(
E_{j}\right)  \approx N^{-1/d} \ . \label{Ej}%
\end{equation}
Let $0<a<b<1/2$, let $f$ be a non-negative function in $L^{2}\left(
\mathbb{T}^{d}\right)  $, and let%
\begin{align}
J\left(  N,f,r\right)   &  =N\int_{E_{1}}\cdots N\int_{E_{N}}\int
_{\mathbb{T}^{d}} \times\left\vert \frac{1}{N}\sum_{j=1}^{N}f\left(
z_{j}\right)  \chi_{-x+B_{r}}\left(  z_{j}\right)  \right. \label{JNf}\\
&  \qquad\left.  -\int_{-x+B_{r}}f\left(  y\right)  dy\right\vert ^{2}%
dxdz_{1}\cdots dz_{N}\nonumber
\end{align}
denote the discrepancy of the so-called jittered sampling. Then there exist
$c_{1}>0$, independent of $N$ and $f$, such that%
\begin{equation}
\int_{a}^{b}J\left(  N,f,r\right)  \ dr\geqslant c_{1}N^{-1-1/d}\left\Vert
f\right\Vert _{2}^{2}\ . \label{ineq_basso}%
\end{equation}
Moreover there exist $c_{2}>0$, independent of $N$, $f$ and $r$, such that
\begin{equation}
J\left(  N,f,r\right)  \leqslant c_{2}N^{-1}\left\Vert f\right\Vert _{2}%
^{2}\ . \label{ineq_alto}%
\end{equation}
These inequalities are best possible. That is, there exist positive constants
$c_{3}$ and $c_{4}$ such that for every $N$ there exist functions $f_{1}$ and
$f_{2}$ in $L^{2}\left(  \mathbb{T}^{d}\right)  $ such that, for every $r$,%
\begin{equation}
J\left(  N,f_{1},r\right)  \leqslant c_{3}N^{-1-1/d}\left\Vert f_{1}%
\right\Vert _{2}^{2}\ ,\ \ \text{\ \ \ }J\left(  N,f_{2},r\right)  \geqslant
c_{4}N^{-1}\left\Vert f_{2}\right\Vert _{2}^{2}\ . \label{best}%
\end{equation}

\end{theorem}

Finally, the upper bound $J\left(  N,f,r\right)  \leqslant c_{2}%
N^{-1}\left\Vert f\right\Vert _{2}^{2}$ in (\ref{ineq_alto}), which holds true
for every function in $L^{2}\left(  \mathbb{T}^{d}\right)  $, can be improved
under suitable regularity assumptions on $f$. Indeed we have the following result.

\begin{theorem}
\label{jitt j}We keep the notation of the previous theorem. Let $f$ be a
non-negative and H\"{o}lder continuous function of order $0<\beta\leqslant1$
on $\mathbb{T}^{d}$, that is $\left\vert f\left(  x\right)  -f\left(
y\right)  \right\vert \leqslant\left\vert x-y\right\vert ^{\beta}$ for all
$x,y\in\mathbb{T}^{d}$. Then, for every $r$,
\[
J\left(  N,f,r\right)  \leqslant\left\{
\begin{array}
[c]{cc}%
cN^{-1-2\beta/d} & \ \text{if }\beta<1/2,\\
cN^{-1-1/d} & \ \text{if }\beta\geqslant1/2.
\end{array}
\right.
\]

\end{theorem}

\section{Proofs}

The starting point in our proofs is the Fourier expansion of the discrepancy
function $x\mapsto\mathcal{D}_{N}\left(  x,r\right)  $. For every
$m\in\mathbb{Z}^{d}$ the Fourier coefficients of this function are given by
\begin{equation}
\widehat{\mathcal{D}}_{N}\left(  m,r\right)  =\int_{\mathbb{T}^{d}}%
\mathcal{D}_{N}\left(  x,r\right)  e^{-2\pi imx}dx \ . \label{not four}%
\end{equation}

%We recall that the function $f$ is real-valued and therefore $\widehat{f}\left(  -m\right) =\overline{f(m)}$ for every $m \in \mathbb{Z}^d$.
%\smallskip

\begin{lemma}
\label{Lemma Parseval}For every $m\in\mathbb{Z}^{d}$ we have%
\begin{equation}
\widehat{\mathcal{D}}_{N}\left(  m,r\right)  =\left(  \frac{1}{N}\sum
_{j=1}^{N}\alpha_{j}e^{2\pi imz_{j}}-\widehat{f}\left(  -m\right)  \right)
\widehat{\chi}_{B_{r}}\left(  m\right)  \ . \label{transf}%
\end{equation}
Therefore, by Parseval's formula,%
\begin{equation}
\int_{\mathbb{T}^{d}}\left\vert \widehat{\mathcal{D}}_{N}\left(  x,r\right)
\right\vert ^{2}dx=\sum_{m\in\mathbb{Z}^{d}}\left\vert \frac{1}{N}\sum
_{j=1}^{N}\alpha_{j}e^{2\pi imz_{j}}-\widehat{f}\left(  -m\right)  \right\vert
^{2}\left\vert \widehat{\chi}_{B_{r}}\left(  m\right)  \right\vert ^{2}\ .
\label{Parseval}%
\end{equation}
Moreover, for every $0<a<b<1/2$,%
\begin{align}
&  \int_{a}^{b}\int_{\mathbb{T}^{d}}\left\vert \widehat{\mathcal{D}}%
_{N}\left(  x,r\right)  \right\vert ^{2}\,dx\,dr\label{Parseval2}\\
&  \approx\sum_{m\in\mathbb{Z}^{d}}\left(  1+\left\vert m\right\vert \right)
^{-d-1}\left\vert \frac{1}{N}\sum_{j=1}^{N}\alpha_{j}e^{2\pi imz_{j}}%
-\widehat{f}\left(  -m\right)  \right\vert ^{2}\ .\nonumber
\end{align}

\end{lemma}

\begin{proof}
By the definition of $\mathcal{D}_{N}\left(  x,r\right)  $ (see Definition
\ref{defdiscr}) we have%
\begin{align*}
&  \widehat{\mathcal{D}}_{N}\left(  m,r\right)  =\int_{\mathbb{T}^{d}%
}\mathcal{D}_{N}\left(  x,r\right)  e^{-2\pi imx}dx\\
&  =\frac{1}{N}\sum_{j=1}^{N}\alpha_{j}\int_{\mathbb{T}^{d}}\chi_{-x+B_{r}%
}\left(  z_{j}\right)  e^{-2\pi imx}dx\\
&  \qquad-\int_{\mathbb{T}^{d}}\int_{\mathbb{T}^{d}}\chi_{-x+B_{r}}\left(
y\right)  f\left(  y\right)  dy\ e^{-2\pi imx}dx\\
&  =\frac{1}{N}\sum_{j=1}^{N}\alpha_{j}\int_{\mathbb{T}^{d}}\chi_{B_{r}%
}\left(  x+z_{j}\right)  e^{-2\pi imx}dx\\
&  \qquad-\int_{\mathbb{T}^{d}}\int_{\mathbb{T}^{d}}\chi_{B_{r}}\left(
x+y\right)  f\left(  y\right)  e^{-2\pi imx}dydx\\
&  =\frac{1}{N}\sum_{j=1}^{N}\alpha_{j}e^{2\pi imz_{j}}\int_{\mathbb{T}^{d}%
}\chi_{B_{r}}\left(  x+z_{j}\right)  e^{-2\pi im\left(  x+z_{j}\right)  }dx\\
&  \qquad-\int_{\mathbb{T}^{d}}f\left(  y\right)  e^{2\pi imy}\int
_{\mathbb{T}^{d}}\chi_{B_{r}}\left(  x+y\right)  e^{-2\pi im\left(
x+y\right)  }dxdy\\
&  =\frac{1}{N}\sum_{j=1}^{N}\alpha_{j}e^{2\pi imz_{j}}\widehat{\chi}_{B_{r}%
}\left(  m\right)  -\widehat{f}\left(  -m\right)  \widehat{\chi}_{B_{r}%
}\left(  m\right)  \ .
\end{align*}

Then Parseval's identity yields (\ref{Parseval}). Since
\begin{equation}
\widehat{\chi}_{B_{r}}\left(  m\right)  =r^{d} \ \frac{J_{d/2}\left(
2\pi\left\vert rm\right\vert \right)  }{\left\vert rm\right\vert ^{d/2}} \ ,
\label{bessel}%
\end{equation}
the asymptotic expansions of the Bessel functions (see e.g. \cite{Guido})
readily give
\begin{equation}
\int_{a}^{b}\left\vert \widehat{\chi}_{B_{r}}\left(  m\right)  \right\vert
^{2}dr\approx\left\vert m\right\vert ^{-d-1} \label{decmed}%
\end{equation}
for every $0\neq m\in\mathbb{Z}^{d}$. Then we obtain (\ref{Parseval2}).
\end{proof}

We also need the following simple inequality. For all complex numbers $a$,
$b$, we have
\begin{equation}
\left\vert a-b\right\vert ^{2}\geqslant\frac{1}{2}\left\vert a\right\vert
^{2}-\left\vert b\right\vert ^{2} \ . \label{oops}%
\end{equation}
The proof of (\ref{oops}) is straightforward,%
\[
\left\vert a\right\vert ^{2}\leqslant\left(  \left\vert a-b\right\vert
+\left\vert b\right\vert \right)  ^{2}=\left\vert a-b\right\vert
^{2}+\left\vert b\right\vert ^{2}+2\left\vert a-b\right\vert \left\vert
b\right\vert \leqslant2\left\vert a-b\right\vert ^{2}+2\left\vert b\right\vert
^{2} \ .
\]

\medskip

\begin{proof}
[Proof of Theorem \ref{five}]Let $k \in L^{1}\left(  \mathbb{R}^{d}\right)  $
be a radial non-negative function such that its Fourier transform
\[
\widehat{k}\left(  \xi\right)  =\int_{\mathbb{R}^{d}}k\left(  x\right)
e^{-2\pi ix\xi}dx
\]
is smooth, with $0\leq\widehat{k}\left(  \xi\right)  \leq1$, with $\widehat
{k}\left(  0\right)  =1$, and $\widehat{k}\left(  \xi\right)  =0$ if
$\left\vert \xi\right\vert \geqslant1$. Let $M$ be a positive integer to be
chosen later, and define the periodic function
\begin{equation}
K_{M}\left(  x\right)  =\sum_{m\in\mathbb{Z}^{d}}M^{d}k\left(  M\left(
x+m\right)  \right)  =\sum_{m\in\mathbb{Z}^{d}}\widehat{k}\left(
M^{-1}m\right)  e^{2\pi imx}. \label{K_M}%
\end{equation}
Observe that
\begin{equation}
K_{M}\left(  0\right)  \geqslant k\left(  0\right)  M^{d}\ . \label{dashhhhhh}%
\end{equation}
Moreover, by the assumptions on $\widehat{k}\left(  \xi\right)  $ we have
$k\left(  0\right)  >0$ and
\[
0\leq\widehat{K}_{M}\left(  m\right)  =\widehat{k}\left(  M^{-1}m\right)
\leq1 \ .
\]
Since we are going to prove an estimate from below, then we may assume $f$
real-valued. Then (\ref{Parseval2}), (\ref{oops}) and (\ref{decmed}) yield%
\begin{align}
&  \int_{a}^{b}\int_{\mathbb{T}^{d}}\left\vert \mathcal{D}_{N}\left(
x,r\right)  \right\vert ^{2}dxdr\label{uno e due}\\
&  \approx\sum_{m\in\mathbb{Z}^{d}}\left(  1+\left\vert m\right\vert \right)
^{-d-1}\left\vert \frac{1}{N}\sum_{j=1}^{N}\alpha_{j}e^{2\pi imz_{j}}%
-\widehat{f}\left(  -m\right)  \right\vert ^{2}\nonumber\\
&  \geqslant c\sum_{m\in\mathbb{Z}^{d}}\widehat{K}_{M}\left(  m\right)
\left(  1+\left\vert m\right\vert \right)  ^{-d-1}\left\vert \frac{1}{N}%
\sum_{j=1}^{N}\alpha_{j}e^{2\pi imz_{j}}-\widehat{f}\left(  -m\right)
\right\vert ^{2}\nonumber\\
&  \geqslant c\left(  1+M\right)  ^{-d-1}\sum_{\left\vert m\right\vert
\leqslant M}\widehat{K}_{M}\left(  m\right)  \left(  \frac{1}{2}\left\vert
\frac{1}{N}\sum_{j=1}^{N}\alpha_{j}e^{2\pi imz_{j}}\right\vert ^{2}-\left\vert
\widehat{f}\left(  m\right)  \right\vert ^{2}\right) \nonumber\\
&  \geqslant cM^{-d-1}\left(  \frac{1}{2}\sum_{m\in\mathbb{Z}^{d}}\widehat
{K}_{M}\left(  m\right)  \left\vert \frac{1}{N}\sum_{j=1}^{N}\alpha_{j}e^{2\pi
imz_{j}}\right\vert ^{2}-\sum_{\left\vert m\right\vert \leqslant M}\left\vert
\widehat{f}\left(  m\right)  \right\vert ^{2}\right) \nonumber\\
&  =cM^{-d-1}\left(  \frac{1}{2}I-II\right)  \ .\nonumber
\end{align}
The non-negativity of the $\alpha_{j}$ and (\ref{dashhhhhh}) yield
\begin{align}
I  &  =\sum_{m\in\mathbb{Z}^{d}}\widehat{K}_{M}\left(  m\right)  \left\vert
\frac{1}{N}\sum_{j=1}^{N}\alpha_{j}e^{2\pi imz_{j}}\right\vert ^{2}%
\label{hugh}\\
&  \geqslant\sum_{m\in\mathbb{Z}^{d}}\widehat{K}_{M}\left(  m\right)  \frac
{1}{N^{2}}\sum_{k=1}^{N}\sum_{j=1}^{N}\alpha_{j}e^{2\pi imz_{j}}\alpha
_{k}e^{-2\pi imz_{k}}\nonumber\\
&  =\frac{1}{N^{2}}\sum_{j,k=1}^{N}\alpha_{j}\alpha_{k}K_{M}\left(
z_{j}-z_{k}\right)  \geqslant\frac{1}{N^{2}}\sum_{k=1}^{N}\alpha_{k}^{2}%
K_{M}\left(  0\right) \nonumber\\
&  \geqslant k\left(  0\right)  M^{d}N^{-1}\left\Vert \alpha\right\Vert
^{2}\ .\nonumber
\end{align}
(The above nice argument is essentially due to H. Montgomery, see
\cite[Theorem 5.12]{montgomery}). In order to deal with $II$ we set $r=q/2$
and $1/r+1/r^{\prime}=1$ (hence $r^{\prime}=q/(q-2)$ and $r^{\prime}=\infty$
when $q=2$). Then the inequalities of H\"{o}lder and Hausdorff--Young yield%
\begin{align*}
II  &  =\sum_{\left\vert m\right\vert \leqslant M}\left\vert \widehat
{f}\left(  m\right)  \right\vert ^{2}\leqslant\left\{  \sum_{\left\vert
m\right\vert \leqslant M}1\right\}  ^{1/r^{\prime}}\left\{  \sum_{\left\vert
m\right\vert \leqslant M}\left\vert \widehat{f}\left(  m\right)  \right\vert
^{2r}\right\}  ^{1/r}\\
&  \leqslant cM^{d/r^{\prime}}\left\{  \sum_{\left\vert m\right\vert \leqslant
M}\left\vert \widehat{f}\left(  m\right)  \right\vert ^{q}\right\}
^{2/q}\leqslant cM^{d\left(  1-2/q\right)  }\left\Vert f\right\Vert _{p}%
^{2}\ .
\end{align*}
Then%
\begin{align*}
\frac{1}{2}I-II  &  \geqslant\frac{1}{2}k\left(  0\right)  M^{d}%
N^{-1}\left\Vert \alpha\right\Vert ^{2}-cM^{d\left(  1-2/q\right)  }\left\Vert
f\right\Vert _{p}^{2}\\
&  =M^{d}\left(  \frac{1}{2}k\left(  0\right)  N^{-1}\left\Vert \alpha
\right\Vert ^{2}-cM^{-2d/q}\left\Vert f\right\Vert _{p}^{2}\right)  .
\end{align*}
Now we set%
\[
M=\left(  \frac{4cN\left\Vert f\right\Vert _{p}^{2}}{k\left(  0\right)
\left\Vert \alpha\right\Vert ^{2}}\right)  ^{q/(2d)},
\]
so that%
\[
cNM^{-2d/q}\left\Vert f\right\Vert _{p}^{2}=\frac{1}{4}k\left(  0\right)
\left\Vert \alpha\right\Vert ^{2}\ .
\]
Then we obtain%
\begin{align*}
&  \int_{a}^{b}\int_{\mathbb{T}^{d}}\left\vert \mathcal{D}_{N}\left(
x,r\right)  \right\vert ^{2}dxdr\geqslant c_{1}M^{-d-1}M^{d}N^{-1}\left\Vert
\alpha\right\Vert ^{2}\\
&  =c_{1}N^{-1}\left(  \frac{4cN\left\Vert f\right\Vert _{p}^{2}}{k\left(
0\right)  \left\Vert \alpha\right\Vert ^{2}}\right)  ^{-q/(2d)}\left\Vert
\alpha\right\Vert ^{2}\\
&  =cN^{-1-q/(2d)}\left\Vert \alpha\right\Vert ^{2+q/d}\left\Vert f\right\Vert
_{p}^{-q/d}.
\end{align*}

\end{proof}

\medskip

\begin{proof}
[Proof of Theorem \ref{Lpalto} ]By Lemma \ref{Lemma Parseval}%
\[
\int_{\mathbb{T}^{d}}\left\vert \widehat{\mathcal{D}}_{N}\left(  x,r\right)
\right\vert ^{2}dx=\sum_{m\in\mathbb{Z}^{d}}\left\vert \frac{1}{N}\sum
_{j=1}^{N}\alpha_{j}e^{2\pi imz_{j}}-\widehat{f}\left(  -m\right)  \right\vert
^{2}\left\vert \widehat{\chi}_{B_{r}}\left(  m\right)  \right\vert ^{2} \ .
\]
We now choose
\[
\left\{  z_{1},\ldots,z_{N}\right\}  =\frac{1}{H}\mathbb{Z}^{d}\cap\left[
-\frac{1}{2},\frac{1}{2}\right)  ^{d}%
\]
and $\alpha_{j}=1$ for every $j$ (hence $\left\Vert \alpha\right\Vert =1$). We
observe that
\[
\sum_{j=1}^{N}e^{2\pi im\cdot z_{j}}=%
\begin{cases}
N & \text{if $m=Hk$ with $k\in\mathbb{Z}^{d}$,}\\
0 & \text{otherwise.}%
\end{cases}
\]
Let $M>0$ and let $F$ be any function in $L^{p}\left(  \mathbb{T}^{d}\right)
$ satisfying $\widehat{F}\left(  m\right)  =1$ if $\left\vert m \right\vert
\leqslant M$. Let $f\left(  x\right)  =F\left(  Hx\right)  $. Since $H$ is a
positive integer, then $\left\Vert f\right\Vert _{p}=\left\Vert F\right\Vert
_{p}$. Also, observe that%
\[
f\left(  x\right)  =F\left(  Hx\right)  =\sum_{k\in\mathbb{Z}^{d}}\widehat
{F}\left(  k\right)  e^{2\pi iHkx} \ .
\]
Hence
\[
\widehat{f}\left(  m\right)  =%
\begin{cases}
\widehat{F}\left(  H^{-1}m\right)  & \text{if $m=Hk$, with $k\in\mathbb{Z}%
^{d}$,}\\
0 & \text{otherwise,}%
\end{cases}
\]
and%
\begin{align*}
&  \sum_{m\in\mathbb{Z}^{d}}\left\vert \frac{1}{N}\sum_{j=1}^{N}e^{2\pi
imz_{j}}-\widehat{f}\left(  -m\right)  \right\vert ^{2}\left\vert
\widehat{\chi}_{B_{r}}\left(  m\right)  \right\vert ^{2}\\
&  =\sum_{m\in\mathbb{Z}^{d}}\left\vert 1-\widehat{F}\left(  -m\right)
\right\vert ^{2}\left\vert \widehat{\chi}_{B_{r}}\left(  Hm\right)
\right\vert ^{2} \ .
\end{align*}
The above equality holds for every $F\left(  x\right)  $. Now we choose a
smooth function $\phi\left(  x\right)  $ on $\mathbb{R}^{d}$ satisfying
\[
\widehat{\phi}\left(  \xi\right)  =%
\begin{cases}
1 & \text{if }\left\vert \xi\right\vert \leqslant1\text{,}\\
0 & \text{if }\left\vert \xi\right\vert >2\text{.}%
\end{cases}
\]
Let
\begin{equation}
F\left(  x\right)  =\sum_{k\in\mathbb{Z}^{d}}M^{d}\phi\left(  M\left(
x+k\right)  \right)  =\sum_{m\in\mathbb{Z}^{d}}\widehat{\phi}\left(
m/M\right)  e^{2\pi imx} \ . \label{dotdot}%
\end{equation}
Then, by the fast decay of $\phi\left(  x\right)  $, i.e. $\left\vert
\phi\left(  x\right)  \right\vert \leqslant c\left(  1+\left\vert x\right\vert
\right)  ^{-L}$ for every $L>0$, where the above constant $c$ depends on $L$,
we have
\begin{align*}
\left\Vert F\right\Vert _{p}  &  \leqslant\sum_{k\in\mathbb{Z}^{d}}\left\Vert
M^{d}\phi\left(  M\left(  x+k\right)  \right)  \right\Vert _{p}\\
&  \leqslant\left\Vert M^{d}\phi\left(  Mx\right)  \right\Vert _{p}%
+cM^{d-L}\sum_{0\neq k\in\mathbb{Z}^{d}}\left\vert k\right\vert ^{-L}\\
&  \leqslant cM^{d}\left(  \int_{\mathbb{R}^{d}}\left\vert \phi\left(
Mx\right)  \right\vert ^{p}dx\right)  ^{1/p}+cM^{d-L}\\
&  \leqslant cM^{d\left(  1-1/p\right)  }\left(  \int_{\mathbb{R}^{d}%
}\left\vert \phi\left(  y\right)  \right\vert ^{p}dy\right)  ^{1/p}%
+cM^{d-L}\leqslant cM^{d/q}\ .
\end{align*}
Since $\widehat{F}\left(  m\right)  =1$ for $\left\vert m\right\vert \leqslant
M$, then (\ref{decmed}) yields
\begin{align*}
\sum_{m\in\mathbb{Z}^{d}}\left\vert 1-\widehat{F}\left(  -m\right)
\right\vert ^{2}\left\vert \widehat{\chi}_{B_{r}}\left(  Hm\right)
\right\vert ^{2}  &  \leqslant c\sum_{\left\vert m\right\vert \geqslant
M}\left\vert \widehat{\chi}_{B_{r}}\left(  Hm\right)  \right\vert ^{2}\\
&  \leqslant cH^{-d-1}M^{-1}\ .
\end{align*}
Since $H^{d}=N$ and since $M^{-1}\leqslant\left\Vert F\right\Vert _{p}^{-q/d}$
one finally obtains%
\[
\sum_{m\in\mathbb{Z}^{d}}\left\vert 1-\widehat{F}\left(  m\right)  \right\vert
^{2}\left\vert \widehat{\chi}_{B_{r}}\left(  Hm\right)  \right\vert
^{2}\leqslant cN^{-1-1/d}\left\Vert f\right\Vert _{p}^{-q/d} \ .
\]
Now we prove the second part of this theorem. For every $N$ let $f$ be a
trigonometric polynomial on $\mathbb{T}^{d}$ such that
\[
\widehat{f}(m)=1\ ,\ \ \ \ \text{if}\ 0\leqslant|m|\leqslant\frac
{1}{\varepsilon\left(  N\right)  }
\]
and
\[
\left\Vert f\right\Vert _{1}\leqslant c\ .
\]
The above function $f$ can be constructed in many ways: we can start with the
previous function $\phi\left(  x\right)  $, or we can take the product of $d$
copies of the $1$-dimensional de la Vall\'{e}e Poussin kernels (see e.g.
\cite[Chapter 1]{Katz04}). For every $j=1,\ldots,N$ let $z_{j}=0$ and
$\alpha_{j}=1$. Then (\ref{bessel}) yields the upper bound
\begin{align*}
&  \int_{\mathbb{T}^{d}}\left\vert \mathcal{D}_{N}\left(  x,r\right)
\right\vert ^{2}dx\leqslant c\sum_{m\in\mathbb{Z}^{d}}|m|^{-\left(
d+1\right)  }\left\vert \frac{1}{N}\sum_{j=1}^{N}\alpha_{j}e^{2\pi imz_{j}%
}-\widehat{f}\left(  -m\right)  \right\vert ^{2}\\
&  \leqslant c\sum_{\left\vert m\right\vert \leqslant1/\varepsilon(N)
}|m|^{-\left(  d+1\right)  }\left\vert 1-1\right\vert ^{2}+c\sum_{\left\vert
m\right\vert >1/\varepsilon(N) }|m|^{-\left(  d+1\right)  }\left\vert
1-\widehat{f}\left(  -m\right)  \right\vert ^{2}\\
&  \leqslant c\sum_{\left\vert m\right\vert >1/\varepsilon(N) }|m|^{-\left(
d+1\right)  }\leqslant c \ \varepsilon(N) \ .
\end{align*}

\end{proof}

\begin{remark}
\label{segni} Let us briefly discuss the assumption $\alpha_{j} \geqslant0$.
For instance, if the $\alpha_{j} $ are not all zero and we assume $\sum
_{j=1}^{N}\alpha_{j}=0$, together with $z_{j}=0$ for every $j$, then
\begin{align*}
&  \int_{a}^{b}\int_{\mathbb{T}^{d}}\left\vert \mathcal{D}_{N}\left(
x,r\right)  \right\vert ^{2}dxdr\\
&  \approx\sum_{m\in\mathbb{Z}^{d}}\left(  1+\left\vert m\right\vert \right)
^{-d-1}\left\vert \frac{1}{N}\sum_{j=1}^{N}\alpha_{j}e^{2\pi imz_{j}}%
-\widehat{f}\left(  -m\right)  \right\vert ^{2}\\
&  =\sum_{m\in\mathbb{Z}^{d}}\left(  1+\left\vert m\right\vert \right)
^{-d-1}\left\vert \widehat{f}\left(  -m\right)  \right\vert ^{2}\ .
\end{align*}
Now let $f\left(  x\right)  =e^{2\pi ikx}$. Then $\left\Vert f\right\Vert
_{p}=1$ and
\[
\int_{a}^{b}\int_{\mathbb{T}^{d}}\left\vert \mathcal{D}_{N}\left(  x,r\right)
\right\vert ^{2}dxdr\approx\left(  1+\left\vert k\right\vert \right)  ^{-d-1}
\ .
\]
This quantity can be made arbitrarily small as $\left\vert k\right\vert $ gets
large.
%Finally, by continuity the same is true if $\sum_{j=1}^{N}%
%\alpha_{j}$ is close to $0$ and each $z_{j}$ is close to zero.

\end{remark}

\medskip

\begin{proof}
[Proof of Theorem \ref{morreybasso}]Let $M>0$ be a constant to be chosen
later. By (\ref{Parseval2}) we have%
\begin{align*}
&  \int_{a}^{b}\int_{\mathbb{T}^{d}}\left\vert \widehat{\mathcal{D}}%
_{N}\left(  x,r\right)  \right\vert ^{2}\,dx\,dr\\
&  \approx\sum_{m\in\mathbb{Z}^{d}}\left(  1+\left\vert m\right\vert \right)
^{-d-1}\left\vert \frac{1}{N}\sum_{j=1}^{N}\alpha_{j}e^{2\pi imz_{j}}%
-\widehat{f}\left(  -m\right)  \right\vert ^{2}\\
&  \geq\left(  1+M\right)  ^{-d-1}\sum_{\left\vert m\right\vert \leq
M}\left\vert \frac{1}{N}\sum_{j=1}^{N}\alpha_{j}e^{2\pi imz_{j}}-\widehat
{f}\left(  -m\right)  \right\vert ^{2}\ .
\end{align*}
Let $K_{M}\left(  x\right)  $ be as in (\ref{K_M}). Observe that, for every
$h>0$, there exists a positive constant $c$ such that%
\begin{equation}
0\leqslant K_{M}\left(  x\right)  \leqslant c \ M^{d} \left(  1+M \left|  x
\right|  \right)  ^{-h} \ . \label{Stime K_M}%
\end{equation}
If $h>\lambda$, then (\ref{Stime K_M}) yields
\begin{align}
&  0\leqslant K_{M}\ast f\left(  x\right) \label{long}\\
&  =\int_{\mathbb{T}^{d}}K_{M}\left(  y\right)  f\left(  x-y\right)
dy\nonumber\\
&  =\int_{\left\{  \left\vert y\right\vert \leqslant M^{-1}\right\}  }%
K_{M}\left(  y\right)  f\left(  x-y\right)  dy\nonumber\\
&  \qquad+\sum_{k=0}^{+\infty}\int_{\left\{  2^{k}M^{-1}<\left\vert
y\right\vert \leqslant2^{k+1}M^{-1}\right\}  \cap\mathbb{T}^{d}}K_{M}\left(
y\right)  f\left(  x-y\right)  dy\nonumber\\
&  \leqslant cM^{d}\int_{\left\{  \left\vert y\right\vert \leqslant
M^{-1}\right\}  }f\left(  x-y\right)  dy\nonumber\\
&  \qquad+cM^{d}\sum_{k=0}^{+\infty}2^{-hk}\int_{\left\{  \left\vert
y\right\vert \leqslant2^{k+1}M^{-1}\right\}  \cap\mathbb{T}^{d}}f\left(
x-y\right)  dy\nonumber\\
&  \leqslant cM^{d-\lambda}\frac{1}{\left(  M^{-1}\right)  ^{\lambda}}%
\int_{\left\{  \left\vert y\right\vert \leqslant M^{-1}\right\}  }\left\vert
f\left(  x-y\right)  \right\vert dy\nonumber\\
&  \qquad+cM^{d}\sum_{k=0}^{+\infty}2^{-hk}\left(  2^{k+1}M^{-1}\right)
^{\lambda}\nonumber\\
&  \qquad\times\frac{1}{\left(  2^{k+1}M^{-1}\right)  ^{\lambda}}%
\int_{\left\{  \left\vert y\right\vert \leqslant2^{k+1}M^{-1}\right\}
\cap\mathbb{T}^{d}}\left\vert f\left(  x-y\right)  \right\vert dy\nonumber\\
&  \leqslant cM^{d-\lambda}\left\Vert f\right\Vert _{1,\lambda}+c2^{\lambda
}M^{d-\lambda}\sum_{k=0}^{+\infty}2^{-\left(  h-\lambda\right)  k}\left\Vert
f\right\Vert _{1,\lambda}\nonumber\\
&  \leqslant cM^{d-\lambda}\left\Vert f\right\Vert _{1,\lambda}\ .\nonumber
\end{align}
Since $0\leqslant\widehat{K}_{M}\left(  m\right)  \leqslant1$, we have%
\begin{align}
&  \sum_{\left\vert m\right\vert \leqslant M}\left\vert \frac{1}{N}\sum
_{j=1}^{N}\alpha_{j}e^{2\pi imz_{j}}-\widehat{f}\left(  -m\right)  \right\vert
^{2}\nonumber\\
&  \geqslant\sum_{\left\vert m\right\vert \leqslant M}\widehat{K}_{M}\left(
m\right)  \left\vert \frac{1}{N}\sum_{j=1}^{N}\alpha_{j}e^{2\pi imz_{j}%
}-\widehat{f}\left(  -m\right)  \right\vert ^{2}\nonumber\\
&  =\sum_{m\in\mathbb{Z}^{d}}\widehat{K}_{M}\left(  m\right)  \left(  \frac
{1}{N}\sum_{j=1}^{N}\alpha_{j}e^{2\pi imz_{j}}-\widehat{f}\left(  -m\right)
\right) \nonumber\\
&  \qquad\times\left(  \frac{1}{N}\sum_{k=1}^{N}\alpha_{k}e^{-2\pi imz_{k}%
}-\overline{\widehat{f}\left(  -m\right)  } \right) \nonumber\\
&  =\frac{1}{N^{2}}\sum_{j,k=1}^{N}\sum_{m\in\mathbb{Z}^{d}}\widehat{K}%
_{M}\left(  m\right)  \alpha_{k}\alpha_{j}e^{2\pi im\left(  z_{j}%
-z_{k}\right)  }\nonumber\\
&  \qquad-\frac{1}{N}\sum_{j=1}^{N}\alpha_{j}\sum_{m\in\mathbb{Z}^{d}}%
\widehat{K}_{M}\left(  m\right)  e^{2\pi imz_{j}}\overline{\widehat{f}\left(
-m\right)  }\nonumber\\
&  \qquad-\frac{1}{N}\sum_{k=1}^{N}\alpha_{k}\sum_{m\in\mathbb{Z}^{d}}%
\widehat{K}_{M}\left(  m\right)  \widehat{f}\left(  -m\right)  e^{-2\pi
imz_{k}}\nonumber\\
&  \qquad+\sum_{m\in\mathbb{Z}^{d}}\widehat{K_{M}}\left(  m\right)  \left\vert
\widehat{f}\left(  m\right)  \right\vert ^{2}\nonumber\\
&  =\frac{1}{N^{2}}\sum_{j,k=1}^{N}\alpha_{k}\alpha_{j}K_{M}\left(
z_{j}-z_{k}\right)  -\frac{2}{N}\mathrm{Re} \left(  \sum_{j=1}^{N}\alpha
_{j}K_{M}\ast f\left(  z_{j}\right)  \right) \label{dd}\\
&  \qquad+\sum_{m\in\mathbb{Z}^{d}}\widehat{K}_{M}\left(  m\right)  \left\vert
\widehat{f}\left(  m\right)  \right\vert ^{2}\ .\nonumber
\end{align}
We are left with the three terms in (\ref{dd}). For the first one, as in
(\ref{hugh}), the non-negativity assumption for the $\alpha_{j}$ yields%
\[
\frac{1}{N^{2}}\sum_{j,k=1}^{N}\alpha_{k}\alpha_{j}K_{M}\left(  z_{j}%
-z_{k}\right)  \geqslant\frac{1}{N^{2}}\sum_{k=1}^{N}\alpha_{k}^{2}%
K_{M}\left(  0\right)  \geqslant c_{1}\frac{1}{N}M^{d}\left\Vert
\alpha\right\Vert ^{2}\ .
\]
As for the second term, (\ref{long}) and the Cauchy-Schwarz inequality yield
\begin{align*}
&  \frac{2}{N}\mathrm{Re} \left(  \sum_{j=1}^{N}\alpha_{j}K_{M}\ast f\left(
z_{j}\right)  \right)  \leqslant\frac{2}{N}\left\vert \sum_{j=1}^{N}\alpha
_{j}K_{M}\ast f\left(  z_{j}\right)  \right\vert \\
&  \leqslant2\left\Vert \alpha\right\Vert \sup_{z \in\mathbb{T}^{d}}\left\vert
K_{M}\ast f\left(  z\right)  \right\vert \leqslant c_{2}\left\Vert
\alpha\right\Vert M^{d-\lambda}\left\Vert f\right\Vert _{1,\lambda}\ .
\end{align*}
We observe that the third term is non-negative. Then
\begin{align*}
&  \sum_{\left\vert m\right\vert \leqslant M}\left\vert \frac{1}{N}\sum
_{j=1}^{N}\alpha_{j}e^{2\pi imz_{j}}-\widehat{f}\left(  -m\right)  \right\vert
^{2}\\
&  \geqslant c_{1}\frac{1}{N}M^{d}\left\Vert \alpha\right\Vert ^{2}%
-c_{2}\left\Vert \alpha\right\Vert M^{d-\lambda}\left\Vert f\right\Vert
_{1,\lambda}\ .
\end{align*}
Setting
\[
M=\left(  \frac{2c_{2}N\left\Vert f\right\Vert _{1,\lambda}}{c_{1}\left\Vert
\alpha\right\Vert }\right)  ^{1/\lambda}%
\]
we have%
\[
c_{1}\frac{1}{N}M^{d}\left\Vert \alpha\right\Vert ^{2}=2c_{2}\left\Vert
\alpha\right\Vert M^{d-\lambda}\left\Vert f\right\Vert _{1,\lambda}\ .
\]
Then
\begin{align}
&  \int_{a}^{b}\int_{\mathbb{T}^{d}}\left\vert \mathcal{D}_{N}\left(
x,r\right)  \right\vert ^{2}dxdr\\
&  \geqslant cM^{-d-1}\sum_{\left\vert m\right\vert \leqslant M}\left\vert
\frac{1}{N}\sum_{j=1}^{N}\alpha_{j}e^{2\pi imz_{j}}-\widehat{f}\left(
-m\right)  \right\vert ^{2}\nonumber\\
&  \geqslant cM^{-1}N^{-1}\left\Vert \alpha\right\Vert ^{2}=cN^{-1-1/\lambda
}\left\Vert f\right\Vert _{1,\lambda}^{-1/\lambda}\left\Vert \alpha\right\Vert
^{2+1/\lambda}\ .\nonumber
\end{align}
Then (\ref{morreysimple}) is proved.
\end{proof}

The proof of Theorem \ref{morreyalto} is not very different from the proof of
Theorem \ref{Lpalto}.

\begin{proof}
[Proof of Theorem \ref{morreyalto}]As in the proof of Theorem \ref{Lpalto},
let $N=H^{d}$, where $H$ is a postive integer. Let
\[
\left\{  z_{1},\ldots,z_{N}\right\}  =\frac{1}{H}\mathbb{Z}^{d}\cap\left[
-\frac{1}{2},\frac{1}{2}\right)  ^{d}\ ,
\]
let $\alpha_{j}=1$ for every $j$, let $M>0$ and let $f\left(  x\right)
=F\left(  Hx\right)  $ with $F \in L^{p}\left(  \mathbb{T}^{d}\right)  $ such
that $\widehat{F}\left(  m\right)  =1$ if $\left\vert m \right\vert \leqslant
M$. As in (\ref{dotdot}) let%
\[
F\left(  x\right)  =\sum_{k\in\mathbb{Z}^{d}}M^{d}\phi\left(  M\left(
x+k\right)  \right)  =\sum_{m\in\mathbb{Z}^{d}}\widehat{\phi}\left(
m/M\right)  e^{2\pi imx}.
\]
By the fast decay of $\phi\left(  x\right)  $, i.e. $\left\vert \phi\left(
x\right)  \right\vert \leqslant c\left(  1+\left\vert x\right\vert \right)
^{-L}$, for every $L>0$, we have
\begin{align*}
\left\Vert F\right\Vert _{1,\lambda}  &  =\sup_{z\in\left[  -1/2,1/2\right)
^{d}, \ 0<r<1/2}\frac{1}{r^{\lambda}}\int_{B\left(  z,r\right)  }\left\vert
M^{d}\sum_{k\in\mathbb{Z}^{d}}\phi\left(  M\left(  x+k\right)  \right)
\right\vert dx\\
&  \leqslant\sup_{z\in\mathbb{R}^{d}, \ r>0}\frac{1}{r^{\lambda}}%
\int_{B\left(  z,r\right)  }\left\vert M^{d}\phi\left(  Mx\right)  \right\vert
dx+cM^{d-L}\sum_{0\neq k\in\mathbb{Z}^{d}}\left\vert k\right\vert ^{-L}\\
&  \leqslant M^{\lambda}\sup_{z\in\mathbb{R}^{d}, \ r>0}\frac{1}{\left(
Mr\right)  ^{\lambda}}\int_{B\left(  z,Mr\right)  }\left\vert \phi\left(
y\right)  \right\vert dy+cM^{d-L}\leqslant cM^{\lambda}\ .
\end{align*}
We claim that%
\[
\left\Vert f\right\Vert _{1,\lambda}\leqslant c\left\Vert F\right\Vert
_{1,\lambda}\leqslant cM^{\lambda}.
\]
Indeed,%
\begin{align*}
\left\Vert f\right\Vert _{1,\lambda}  &  =\sup_{z\in\mathbb{T}^{d},
\ 0<r<1/2}\frac{1}{r^{\lambda}}\int_{B\left(  z,r\right)  }\left\vert F\left(
Hx\right)  \right\vert dx\\
&  =\sup_{z\in\left[  -1/2,1/2\right)  ^{d}, \ 0<r<1/2}\frac{1}{r^{\lambda}%
}\int_{B\left(  Hz,Hr\right)  }\left\vert F\left(  y\right)  \right\vert
H^{-d}dy\\
&  =\sup_{z\in\left[  -1/2,1/2\right)  ^{d}, \ 0<r<1/2}\frac{1}{r^{\lambda}%
}\int_{B\left(  z,Hr\right)  }\left\vert F\left(  y\right)  \right\vert
H^{-d}dy\\
&  =H^{\lambda-d}\sup_{z\in\left[  -1/2,1/2\right)  ^{d}, \ 0<r<1/2}\frac
{1}{\left(  Hr\right)  ^{\lambda}}\int_{B\left(  z,Hr\right)  }\left\vert
F\left(  y\right)  \right\vert dy\ .
\end{align*}
Let $Hr<1/2$, then%
\begin{align*}
&  \frac{1}{\left(  Hr\right)  ^{\lambda}}\int_{B\left(  z,Hr\right)
}\left\vert F\left(  y\right)  \right\vert dy\\
&  \leqslant\sup_{z\in\left[  -1/2,1/2\right)  ^{d}, \ 0<s<1/2}\frac
{1}{s^{\lambda}}\int_{B\left(  z,s\right)  }\left\vert F\left(  y\right)
\right\vert dy\leqslant\left\Vert F\right\Vert _{1,\lambda}.
\end{align*}
Let $Hr\geqslant1/2$. Since $r<1/2$ and $\lambda\leqslant d$, then
\begin{align*}
\frac{1}{\left(  Hr\right)  ^{\lambda}}\int_{B\left(  z,Hr\right)  }\left\vert
F\left(  y\right)  \right\vert dy  &  \leqslant\frac{c}{\left(  Hr\right)
^{\lambda}}\left(  Hr\right)  ^{d}\int_{\left[  -1/2,1/2\right)  ^{d}%
}\left\vert F\left(  y\right)  \right\vert dy\\
&  \leqslant cr^{d-\lambda}H^{d-\lambda}\left\Vert F\right\Vert _{1,\lambda
}\leqslant cH^{d-\lambda}\left\Vert F\right\Vert _{1,\lambda}\ .
\end{align*}
Hence,%
\[
\left\Vert f\right\Vert _{1,\lambda}\leqslant c\max\left(  H^{\lambda
-d}\left\Vert F\right\Vert _{1,\lambda},\left\Vert F\right\Vert _{1,\lambda
}\right)  \leqslant c\left\Vert F\right\Vert _{1,\lambda}\leqslant
cM^{\lambda}\ .
\]
Since $\widehat{F}\left(  m\right)  =1$ for $\left\vert m\right\vert \leqslant
M$, then%
\begin{align*}
\sum_{m\in\mathbb{Z}^{d}}\left\vert 1-\widehat{F}\left(  -m\right)
\right\vert ^{2}\left\vert \widehat{\chi}_{B_{r}}\left(  Hm\right)
\right\vert ^{2}  &  \leqslant c\sum_{\left\vert m\right\vert >M}\left\vert
\widehat{\chi}_{B_{r}}\left(  Hm\right)  \right\vert ^{2}\\
&  \leqslant cH^{-d-1}M^{-1}\ .
\end{align*}
Since $H^{d}=N$ and since $M^{-1}\leqslant c\left\Vert f\right\Vert
_{1,\lambda}^{-1/\lambda}$ one finally obtains%
\[
\sum_{m\in\mathbb{Z}^{d}}\left\vert 1-\widehat{F}\left(  -m\right)
\right\vert ^{2}\left\vert \widehat{\chi}_{B_{r}}\left(  Hm\right)
\right\vert ^{2}\leqslant cN^{-1-1/d}\left\Vert f\right\Vert _{1,\lambda
}^{-1/\lambda}\ .
\]

\end{proof}

\begin{proof}
[Proof of Theorem \ref{jitt}]We recall (see (\ref{transf})) that
\[
\widehat{\mathcal{D}}_{N}\left(  m\right)  =\widehat{\chi}_{B_{r}}\left(
m\right)  \left[  \frac{1}{N}\sum_{j=1}^{N}f\left(  z_{j}\right)  e^{2\pi
im\cdot z_{j}}-\widehat{f}\left(  -m\right)  \right]  \ .
\]
Since $f$ is real-valued we have $\widehat{f}\left(  -m\right)  =\overline
{\widehat{f}\left(  m\right)  }$. Then, see (\ref{JNf}), Parseval's identity
yields
\begin{align}
&  J\left(  N,f,r\right) \nonumber\\
&  =N\int_{E_{1}}\cdots N\int_{E_{N}}\sum_{m\in\mathbb{Z}^{d}}\left\vert
\widehat{\chi}_{B_{r}}\left(  m\right)  \right\vert ^{2}\nonumber\\
&  \qquad\times\left\vert \frac{1}{N}\sum_{j=1}^{N}f\left(  z_{j}\right)
e^{2\pi imz_{j}}-\widehat{f}\left(  -m\right)  \right\vert ^{2}dz_{1}\cdots
dz_{N}\nonumber\\
&  =\sum_{m\in\mathbb{Z}^{d}}\left\vert \widehat{\chi}_{B_{r}}\left(
m\right)  \right\vert ^{2}N\int_{E_{1}}\cdots N\int_{E_{N}}\nonumber\\
&  \qquad\times\left(  \frac{1}{N^{2}}\sum_{j=1}^{N}f^{2}\left(  z_{j}\right)
+\frac{1}{N^{2}}\sum_{j\neq k}f\left(  z_{j}\right)  f\left(  z_{k}\right)
e^{-2\pi im\left(  z_{j}-z_{k}\right)  }\right. \nonumber\\
&  \qquad-\left.  \frac{2}{N}\operatorname{Re}\left(  \overline{\widehat
{f}\left(  -m\right)  }\sum_{j=1}^{N}f\left(  z_{j}\right)  e^{2\pi imz_{j}%
}\right)  +\left\vert \widehat{f}\left(  m\right)  \right\vert ^{2}\right)
dz_{1}\cdots dz_{N}\nonumber\\
&  =\frac{1}{N}\sum_{m\in\mathbb{Z}^{d}}\left\vert \widehat{\chi}_{B_{r}%
}\left(  m\right)  \right\vert ^{2}\sum_{j=1}^{N}\int_{E_{j}}f^{2}\left(
z_{j}\right)  dz_{j}\nonumber\\
&  \qquad+\sum_{m\in\mathbb{Z}^{d}}\left\vert \widehat{\chi}_{B_{r}}\left(
m\right)  \right\vert ^{2}\sum_{j\neq k}\int_{E_{j}}f\left(  z_{j}\right)
e^{-2\pi imz_{j}}dz_{j}\nonumber\\
&  \qquad\times\int_{E_{k}}f\left(  z_{k}\right)  e^{2\pi imz_{k}}%
dz_{k}\nonumber\\
&  \qquad-2\sum_{m\in\mathbb{Z}^{d}}\left\vert \widehat{\chi}_{B_{r}}\left(
m\right)  \right\vert ^{2}\operatorname{Re}\left(  \widehat{f}\left(
m\right)  \sum_{j=1}^{N}\int_{E_{j}}f\left(  z_{j}\right)  e^{2\pi imz_{j}%
}dz_{j}\right) \nonumber\\
&  \qquad+\sum_{m\in\mathbb{Z}^{d}}\left\vert \widehat{\chi}_{B_{r}}\left(
m\right)  \right\vert ^{2}\left\vert \widehat{f}\left(  m\right)  \right\vert
^{2}\nonumber\\
&  =\frac{1}{N}\sum_{m\in\mathbb{Z}^{d}}\left\vert \widehat{\chi}_{B_{r}%
}\left(  m\right)  \right\vert ^{2}\left\Vert f\right\Vert _{2}^{2}+\sum
_{m\in\mathbb{Z}^{d}}\left\vert \widehat{\chi}_{B_{r}}\left(  m\right)
\right\vert ^{2}\sum_{j,k=1}^{N}\left(  \int_{E_{j}}f\left(  z_{j}\right)
e^{-2\pi imz_{j}}dz_{j}\right. \nonumber\\
&  \qquad\left.  \times\int_{E_{k}}f\left(  z_{k}\right)  e^{2\pi imz_{k}%
}dz_{k}\right)  -\sum_{m\in\mathbb{Z}^{d}}\left\vert \widehat{\chi}_{B_{r}%
}\left(  m\right)  \right\vert ^{2}\sum_{j=1}^{N}\left\vert \int_{E_{j}%
}f\left(  z_{j}\right)  e^{-2\pi imz_{j}}dz_{j}\right\vert ^{2}\nonumber\\
&  \qquad-2\sum_{m\in\mathbb{Z}^{d}}\left\vert \widehat{\chi}_{B_{r}}\left(
m\right)  \right\vert ^{2}\left\vert \widehat{f}\left(  m\right)  \right\vert
^{2}+\sum_{m\in\mathbb{Z}^{d}}\left\vert \widehat{\chi}_{B_{r}}\left(
m\right)  \right\vert ^{2}\left\vert \widehat{f}\left(  m\right)  \right\vert
^{2}\nonumber\\
&  =\frac{1}{N}\sum_{m\in\mathbb{Z}^{d}}\left\vert \widehat{\chi}_{B_{r}%
}\left(  m\right)  \right\vert ^{2}\left\Vert f\right\Vert _{2}^{2}+\sum
_{m\in\mathbb{Z}^{d}}\left\vert \widehat{\chi}_{B_{r}}\left(  m\right)
\right\vert ^{2}\left\vert \widehat{f}\left(  m\right)  \right\vert
^{2}\nonumber\\
&  \qquad-\sum_{j=1}^{N}\sum_{m\in\mathbb{Z}^{d}}\left\vert \widehat{\chi
}_{B_{r}}\left(  m\right)  \right\vert ^{2}\left\vert \int_{E_{j}}f\left(
z\right)  e^{-2\pi imz}dz\right\vert ^{2}\nonumber\\
&  \qquad-\sum_{m\in\mathbb{Z}^{d}}\left\vert \widehat{\chi}_{B_{r}}\left(
m\right)  \right\vert ^{2}\left\vert \widehat{f}\left(  m\right)  \right\vert
^{2}\nonumber\\
&  =\frac{1}{N}\left\vert B_{r}\right\vert \left\Vert f\right\Vert _{2}%
^{2}-\sum_{j=1}^{N}\left\Vert \chi_{B_{r}}\ast\left(  f\chi_{E_{j}}\right)
\right\Vert _{2}^{2}\ .\nonumber
\end{align}
In short, we have the nice identity%
\begin{equation}
J\left(  N,f,r\right)  =N^{-1}\left\vert B_{r}\right\vert \left\Vert
f\right\Vert _{2}^{2}-\sum_{j=1}^{N}\left\Vert \chi_{B_{r}}\ast\left(
f\chi_{E_{j}}\right)  \right\Vert _{2}^{2}\ , \label{identitona}%
\end{equation}
which holds true for every $f\in L^{2}\left(  \mathbb{T}^{d}\right)  $. Then
we readily obtain
\[
J\left(  N,f,r\right)  \leqslant cN^{-1}\left\Vert f\right\Vert _{2}^{2}\ .
\]
In order to prove the inequality (\ref{ineq_basso}), that is
\[
\int_{a}^{b}J\left(  N,f,r\right)  \ dr\geqslant c_{1}N^{-1-1/d}\left\Vert
f\right\Vert _{2}^{2}\ ,
\]
we appeal to the lower estimate in (\ref{caso P 2}), which holds true for
\underline{every} choice of the points $z_{1},\ldots,z_{N}$, and then we
integrate $z_{1},\ldots,z_{N}$. Let%
\begin{equation}
f_{z}=\left(  f\left(  z_{1}\right)  ,\ldots,f\left(  z_{N}\right)  \right)
\ . \label{fz}%
\end{equation}
Then the H\"{o}lder inequality yields%
\begin{align*}
&  \int_{a}^{b}J\left(  N,f,r\right)  \ dr\\
&  =N\int_{E_{1}}\cdots N\int_{E_{N}}\int_{\mathbb{T}^{d}}\int_{a}^{b}\\
&  \qquad\times\left\vert \frac{1}{N}\sum_{j=1}^{N}f\left(  z_{j}\right)
\chi_{-x+B_{r}}\left(  z_{j}\right)  -\int_{-x+B_{r}}f\left(  y\right)
dy\right\vert ^{2}drdxdz_{1}\cdots dz_{N}\\
&  \geqslant cN\int_{E_{1}}\cdots N\int_{E_{N}}N^{-1-1/d}\left\Vert
f_{z}\right\Vert ^{2+2/d}\left\Vert f\right\Vert _{2}^{-2/d}dz_{1}\cdots
dz_{N}\\
&  =cN\int_{E_{1}}\cdots N\int_{E_{N}}N^{-1-1/d}\left(  \frac{1}{N}\sum
_{j=1}^{N}f^{2}\left(  z_{j}\right)  \right)  ^{1+1/d}\\
&  \qquad\times\left\Vert f\right\Vert _{2}^{-2/d}dz_{1}\cdots dz_{N}\\
&  \geqslant cN^{-1-1/d}\left\Vert f\right\Vert _{2}^{-2/d}\\
&  \qquad\times\left(  N\int_{E_{1}}\cdots N\int_{E_{N}}\frac{1}{N}\sum
_{j=1}^{N}f^{2}\left(  z_{j}\right)  dz_{1}\cdots dz_{N}\right)  ^{1+1/d}\\
&  =cN^{-1-1/d}\left\Vert f\right\Vert _{2}^{-2/d}\left(  \frac{1}{N}%
\sum_{j=1}^{N}N\int_{E_{j}}f^{2}\left(  z_{j}\right)  dz_{j}\right)
^{1+1/d}\\
&  =cN^{-1-1/d}\left\Vert f\right\Vert _{2}^{-2/d}\left\Vert f\right\Vert
_{2}^{2+2/d}=cN^{-1-1/d}\left\Vert f\right\Vert _{2}^{2}\ .
\end{align*}
So far we have proved (\ref{ineq_basso}) and (\ref{ineq_alto}). We still have
to prove (\ref{best}). We recall that%
\[
J\left(  N,f,r\right)  =N^{-1}\left\vert B_{r}\right\vert \left\Vert
f\right\Vert _{2}^{2}-\sum_{j=1}^{N}\left\Vert \chi_{B_{r}}\ast\left(
f\chi_{E_{j}}\right)  \right\Vert _{2}^{2} \ .
\]
We choose $f$ supported in $E_{1}$ and such that $\left\Vert f\right\Vert
_{2}^{2}\geqslant2N\left\Vert f\right\Vert _{1}^{2}$. Then%
\begin{align*}
&  J\left(  N,f,r\right)  =N^{-1}\left\vert B_{r}\right\vert \left\Vert
f\right\Vert _{2}^{2}-\left\Vert \chi_{B_{r}}\ast f\right\Vert _{2}^{2}\\
&  \geqslant N^{-1}\left\vert B_{r}\right\vert \left\Vert f\right\Vert
_{2}^{2}-\left\Vert \chi_{B_{r}}\right\Vert _{2}^{2}\left\Vert f\right\Vert
_{1}^{2}=\ N^{-1}\left\vert B_{r}\right\vert \left(  \left\Vert f\right\Vert
_{2}^{2}-N\left\Vert f\right\Vert _{1}^{2}\right) \\
&  \geqslant\frac{1}{2}N^{-1}\left\vert B_{r}\right\vert \left\Vert
f\right\Vert _{2}^{2}\ .
\end{align*}
Finally, the existence of a function $f$ satisfying
\[
J\left(  N,f,r\right)  \leqslant c_{3}N^{1-1/d}\left\Vert f\right\Vert
_{2}^{2}
\]
is a particular case of Theorem \ref{jitt j}.
\end{proof}

\medskip

\begin{proof}
[Proof of Theorem \ref{jitt j}]In this proof the integration in $r$ plays no
role. Hence we may forget $r$ and replace $B_{r}$ with $B$. The RHS in (\ref{identitona}) equals%
\[
\sum_{j=1}^{N}\left(  N^{-1}\left\vert B\right\vert \left\Vert f\chi_{E_{j}%
}\right\Vert _{2}^{2}-\left\Vert \chi_{B}\ast\left(  f\chi_{E_{j}}\right)
\right\Vert _{2}^{2}\right)
\]
We are going to show that all the terms in the above sum give the same
contribution.
 For every function
$g\in L^{1}\left(  \mathbb{T}^{d}\right)  $ and every $j=1,\ldots,N$ let%
\begin{equation}
\mathcal{M}_{E_{j}}\left(  g\right)  =\frac{1}{\left\vert E_{j}\right\vert
}\int_{E_{j}}g\ .\label{MQg}%
\end{equation}
 Then, for every $j=1,\ldots,N$%
\begin{align*}
&  N^{-1}\left\vert B\right\vert \left\Vert f\chi_{E_{j}}\right\Vert _{2}%
^{2}-\left\Vert \chi_{B}\ast\left(  f\chi_{E_{j}}\right)  \right\Vert _{2}%
^{2}\\
&  =N^{-2}\left\vert B\right\vert \left(  \mathcal{M}_{E_{j}}\left(
f^{2}\right)  -\mathcal{M}_{E_{j}}^{2}\left(  f\right)  \right)  \\
&  \qquad+N^{-2}\left(  \left\vert B\right\vert \mathcal{M}_{E_{j}}^{2}\left(
f\right)  -\left\Vert \chi_{B}\ast\frac{f\chi_{E_{j}}}{\left\vert
E_{j}\right\vert }\right\Vert _{2}^{2}\right)  \\
&  =N^{-2}\left\vert B\right\vert \left(  \mathcal{M}_{E_{j}}\left(  \left(
f-\mathcal{M}_{E_{j}}\left(  f\right)  \right)  ^{2}\right)  \right.  \\
&  \qquad+N^{-2}\left(  \left\vert B\right\vert \mathcal{M}_{E_{j}}^{2}\left(
f\right)  -\left\Vert \chi_{B}\ast\frac{f\chi_{E_{j}}}{\left\vert
E_{j}\right\vert }\right\Vert _{2}^{2}\right)  \ .
\end{align*}
We recall that $\operatorname*{diam}\left(  E_{j}\right)  \approx N^{-1/d}$.
Since $f$ is H\"{o}lder continuous with exponent $\beta$, then, for any $x\in
E_{j}$, we have%
\[
\left\vert f\left(  x\right)  -\mathcal{M}_{E_{j}}\left(  f\right)
\right\vert \leqslant c\operatorname*{diam}\left(  E_{j}\right)  ^{\beta
}\leqslant cN^{-\beta/d}\ .
\]
Hence
\begin{equation}
N^{-2}\left\vert B\right\vert \mathcal{M}_{E_{j}}\left(  \left(
f-\mathcal{M}_{E_{j}}\left(  f\right)  \right)  ^{2}\right)  \leqslant
cN^{-2-2\beta/d}.\label{Stima primo pezzo}%
\end{equation}
Let us consider the second term, that is
\[
N^{-2}\left(  \left\vert B\right\vert \mathcal{M}_{E_{j}}^{2}\left(  f\right)
-\left\Vert \chi_{B}\ast\frac{f\chi_{E_{j}}}{\left\vert E_{j}\right\vert
}\right\Vert _{2}^{2}\right)  .
\]
Let
$
t_{j}\in E_{j}$ and $B_{j}=t_{j}+B
$.
First observe that if $x\in B_{j}$ and $\operatorname*{dist}\left(  x,\partial
B_{j}\right)  \geqslant cN^{-1/d}$ then
\[
\chi_{B}\ast\frac{f\chi_{E_{j}}}{\left\vert E_{j}\right\vert }\left(
x\right)  =\frac{1}{\left\vert E_{j}\right\vert }\int_{E_{j}}f\left(
y\right)  \ dy=\mathcal{M}_{E_{j}}\left(  f\right)  .
\]
Moreover if $x\notin B_{j}$ and $\operatorname*{dist}\left(  x,\partial
B_{j}\right)  \geqslant cN^{-1/d}$, then%
\[
\chi_{B}\ast\frac{f\chi_{E_{j}}}{\left\vert E_{j}\right\vert }\left(
x\right)  =0\ .
\]
In short, the assumption $\operatorname*{dist}\left(  x,\partial B_{j}\right)
\geqslant cN^{-1/d}$ yields
\begin{equation}
\chi_{B}\ast\frac{f\chi_{E_{j}}}{\left\vert E_{j}\right\vert }\left(
x\right)  =\mathcal{M}_{E_{j}}\left(  f\right)  \chi_{B_{j}}\left(  x\right)
\ .\label{Stima conv}%
\end{equation}
Now let
\[
\mathcal{C}_{j}=\left\{  x\in\mathbb{T}^{d}:\operatorname*{dist}\left(
x,\partial B_{j}\right)  \leqslant cN^{-1/d}\right\}  \ .
\]
Observe that $\left\vert \mathcal{C}_j \right\vert \leqslant cN^{-1/d}$.  Then (\ref{Stima conv}) and the Cauchy-Schwartz inequality \ $\mathcal{M}_{E_{j}}^{2}\left(
f\right)   \leqslant 
\mathcal{M}_{E_{j}}\left(  f^{2}\right)$ \  yield%
\begin{align}
&  \left\vert B\right\vert \mathcal{M}_{E_{j}}^{2}\left(  f\right)
-\left\Vert \chi_{B}\ast\frac{f\chi_{E_{j}}}{\left\vert E_{j}\right\vert
}\right\Vert _{2}^{2}\label{Stima secondo pezzo}\\
&  =\left\vert B\right\vert \mathcal{M}_{E_{j}}^{2}\left(  f\right)  -\left(
\mathcal{M}_{E_{j}}^{2}\left(  f\right)  \left\vert B_j\setminus
\mathcal{C}_{j}\right\vert +\int_{\mathcal{C}_{j}}\left(  \chi_{B}\ast
\frac{f\chi_{E_{j}}}{\left\vert E_{j}\right\vert }\left(  x\right)  \right)
^{2}dx\right)  \nonumber\\
&  =\mathcal{M}_{E_{j}}^{2}\left(  f\right)  \left\vert B_{j}\cap
\mathcal{C}_{j}\right\vert -\int_{\mathcal{C}_{j}}\left(  \chi_{B}\ast
\frac{f\chi_{E_{j}}}{\left\vert E_{j}\right\vert }\left(  x\right)  \right)
^{2}dx\nonumber\\
&  \leqslant\mathcal{M}_{E_{j}}^{2}\left(  f\right)  \left\vert B_{j}%
\cap\mathcal{C}_{j}\right\vert \leqslant\mathcal{M}_{E_{j}}^{2}\left(
f\right)  \left\vert \mathcal{C}_{j}\right\vert \leqslant cN^{-1/d}%
\mathcal{M}_{E_{j}}\left(  f^{2}\right)  \ .\nonumber
\end{align}
Then, using (\ref{Stima primo pezzo}), (\ref{Stima secondo pezzo}), and
summing over all $E_{j}$ we obtain%
\begin{align*}
&  N\int_{E_{1}}\cdots N\int_{E_{N}}\int_{\mathbb{T}^{d}}\left\vert
\mathcal{D}_{N}\left(  x\right)  \right\vert ^{2}dxdz_{1}\cdots dz_{N}\\
&  \leqslant c \sum_{j=1}^{N}\left(  N^{-2-2\beta/d}+N^{-2-1/d}\mathcal{M}%
_{E_{j}}\left(  f^{2}\right)  \right)  \\
&  \leqslant c\left(N^{-1-2\beta/d}+N^{-1-1/d}\left\Vert f\right\Vert _{2}^{2} \right) \ .
\end{align*}
This ends the proof of Theorem \ref{jitt j}.
\end{proof}

\end{document}